\documentstyle[12pt]{article}
\begin{document}\begin{flushright}\thispagestyle{empty}
OUT--4102--80\\
20 June 1999
                                                 \end{flushright}\vspace*{2mm}
                                                 \begin{center}{ \Large\bf
A seventeenth-order polylogarithm ladder$^{a)}$
                                                 }\vglue 10mm{\large\bf
David H.~Bailey$^{b)}$
and David J.~Broadhurst$^{c)}$                   }\end{center}\vfill
                                                 \noindent{\bf Abstract}\quad
Cohen, Lewin and Zagier found four ladders that entail the polylogarithms
${\rm Li}_n(\alpha_1^{-k}):=\sum_{r>0}\alpha_1^{-k r}/r^n$ at order $n=16$,
with indices $k\le360$, and $\alpha_1$ being the smallest known Salem number,
i.e.\ the larger real root of Lehmer's celebrated polynomial
$\alpha^{10}+\alpha^9-\alpha^7-\alpha^6-\alpha^5-\alpha^4-\alpha^3+\alpha+1$,
with the smallest known non-trivial Mahler measure. By adjoining the index
$k=630$, we generate a fifth ladder at order 16 and a ladder at order 17 that
we presume to be unique. This empirical integer relation, between elements of
$\{{\rm Li}_{17}(\alpha_1^{-k})\mid0\le k\le630\}$ and
$\{\pi^{2j}(\log\alpha_1)^{17-2j}\mid 0\le j\le8\}$, entails 125 constants,
multiplied by integers with nearly 300 digits. It has been checked to more
than 59,000 decimal digits. Among the ladders that we found in other number
fields, the longest has order 13 and index 294. It is based on
$\alpha^{10}-\alpha^6-\alpha^5-\alpha^4+1$, which gives the sole Salem
number $\alpha<1.3$ with degree $d<12$ for which $\alpha^{1/2}+\alpha^{-1/2}$
fails to be the largest eigenvalue of the adjacency matrix of a graph.
\vfill\footnoterule\noindent{\small
$^a$) This work was supported by the Director, Office of Computational
and Technology Research, Division of Mathematical, Information, and
Computational Sciences of the U.S. Department of Energy, under contract
number DE-AC03-76SF00098.\\
$^b$) Lawrence Berkeley Laboratory, MS 50B-2239, Berkeley, CA 94720, USA\\
{\tt dhbailey@lbl.gov}\\
$^c)$ Open University, Department of Physics, Milton Keynes MK7 6AA, UK\\
{\tt D.Broadhurst@open.ac.uk}}
\newpage

\section{Introduction}

The findings reported here arose from the discovery by the second author that
\begin{equation}
\alpha^{630}-1={
(\alpha^{315}-1)
(\alpha^{210}-1)
(\alpha^{126}-1)^2
(\alpha^{90}-1)
(\alpha^3-1)^3
(\alpha^2-1)^5
(\alpha-1)^3\over
(\alpha^{35}-1)
(\alpha^{15}-1)^2
(\alpha^{14}-1)^2
(\alpha^5-1)^6
\alpha^{68}}
\label{a630}
\end{equation}
where $\alpha$ is one of the 10 algebraic integers
of Lehmer's remarkable number field~\cite{DHL}
\begin{equation}
\alpha^{10}+\alpha^9-\alpha^7-\alpha^6-\alpha^5-\alpha^4-\alpha^3+\alpha+1
=0
\label{a10}
\end{equation}
Once found, the cyclotomic relation~(\ref{a630}) was proven by
(oft) repeated substitution for $\alpha^{10}$.
It led us to believe that a valid ladder of polylogarithms
exists at order $n=17$, contrary to a suggestion in~\cite{DZ}.
Indeed, we were able to adjoin the index $k=630$ to those with $k\le360$,
found by Henri Cohen, Leonard Lewin and Don Zagier~\cite{CLZ}, and obtain
\begin{equation}
N(n)=77-\lfloor9n/2\rfloor
\label{Nn}
\end{equation}
ladders at orders $n=2\ldots17$.
In particular, at $n=17$,
we found 125 non-zero integers\footnote{See
{\tt ftp://physics.open.ac.uk/pub/physics/dbroadhu/lehmer/integers.txt}}
$a$, $b_j$, $c_k$, with less than
300 digits, such that an empirical relation
\begin{equation}
a\,\zeta(17)=
\sum_{j=0}^8b_j\,\pi^{2j}(\log\alpha_1)^{17-2j}
+\sum_{k\in D({\cal S})}c_k\,{\rm Li}_{17}(\alpha_1^{-k})
\label{result}
\end{equation}
holds to more than 59,000 decimal digits, where
\begin{equation}
\alpha_1=1.176280818259917506544070338474035050693415806564\ldots
\label{alpha}
\end{equation}
is the larger real root of~(\ref{a10}), and the 115 indices $k$
in ${\rm Li}_n(\alpha_1^{-k}):=\sum_{r>0}\alpha_1^{-k r}/r^n$
are drawn from the set, $D({\cal S})$, of positive
integers that divide at least one element of
\begin{eqnarray}
&&{\cal S}:=
\{29,47,50,52,56,57,64,74,75,76,78,84,86,92,96,98,108,110,118,124,130,
\nonumber\\&&
132,138,144,154,160,165,175,182,186,195,204,212,240,246,270,286,360,630\}
\label{DS}
\end{eqnarray}
The coefficient of $\zeta(17)$ was partially factorized as follows
\begin{eqnarray}
a&=&2^7\times3^7\times5^4\times7\times11\times13\times17
\times722063\times15121339
\times379780242109750106753\nonumber\\&&{}
\times5724771750303829791195961
\times C_{217}
\label{ais}
\end{eqnarray}
where
\begin{eqnarray}
&&C_{217}:=5203751052922114540188667952627280712081039342696719260003747081
\nonumber\\&&
41977100981249686783730105404186042839389917052601102889831046723208680
\nonumber\\&&
07066945997308654073833814804516883406394532403532415753816146816138731
\nonumber\\&&
90080853089
\label{C217}
\end{eqnarray}
is a 217-digit non-prime, whose factorization
has not yet been obtained.
The integers in~(\ref{result}) were obtained using less than 4,000 digits
of working precision. The chance of a
numerical accident is thus less than $10^{-55000}$.

In section~2, we review the algorithm for generating such ladders
from cyclotomic relations~\cite{DZ}.
In section~3, we describe our computational strategy, based on the PSLQ
integer relation finder~\cite{ppslq}.
In section~4, we study ladders based on Salem numbers larger
than~(\ref{alpha}), commenting on a connection to graph theory,
observed with Gert Almkvist.

\section{Ladder building in self-reciprocal number fields}

Consider the cyclotomic polynomials, $\Phi_k(x)$,
defined recursively by
\begin{equation}
x^k-1=\prod_{j|k}\Phi_{j}(x)
\label{cycp}
\end{equation}
A real algebraic number $x>1$ is said to satisfy a cyclotomic relation
with index $k$ if there exist rational numbers $\{A_j\mid0\le j<k\}$
such that
\begin{equation}
\Phi_k(x)=x^{A_0}\prod_{j=1}^{k-1}\Phi_{j}^{A_j}(x)
\label{cycr}
\end{equation}
For example, identity~(\ref{a630}) establishes that
$\alpha_1$ -- the smallest known Salem number, and also the smallest known
non-trivial Mahler measure --
satisfies a cyclotomic relation with index $k=630$, which is
63 times larger than the degree of $\alpha_1$. This ratio is larger than
any heretofore discovered, the previous~\cite{CLZ} record being $k/d=36$.

A Salem number is a real algebraic integer, $\alpha>1$, of degree $d=2+2s$,
with $2s>0$ conjugates on the unit circle and the remaining conjugate,
$1/\alpha$, inside it. For any monic polynomial, $P$, with integer
coefficients, the Mahler measure, $M(P)$, is the product of the
absolute values of the roots outside the unit circle.
Thus a Salem number is the Mahler measure of its minimal polynomial.
Derrick Lehmer~\cite{DHL} conjectured that there exists a constant $c>1$
such that $M(P)\ge c$ for all $M(P)>1$.
A stronger form of this conjecture is that $c=\alpha_1$,
the Mahler measure of~(\ref{a10}), found by Lehmer
more than 60 years ago, and still~\cite{DWB1,DWB2,MJM}
the smallest known $M(P)>1$.

In describing how to build a polylogarithmic ladder,
we shall restrict attention to an algebraic number field whose defining
polynomial, $P(x)=x^dP(1/x)$, of even degree $d$, is reciprocal
(i.e.\ palindromic), as in~(\ref{a10}).
Moreover, we require at least one real root with $x>1$. Hence
the discussion encompasses all Salem numbers.
We define cyclotomic norms~\cite{DZ}
\begin{equation}
N_k:=\prod_{r=1}^d\Phi_k(x_r)
\label{norm}
\end{equation}
where the product is over all roots of $P(x)$,
so that $N_k$ is an integer.
Then a necessary condition for a cyclotomic relation of index $k$
is that every prime factor of $N_k$ is also a factor of a norm $N_j$
with $j<k$. This simplifies, considerably,
the business of finding all indices of cyclotomic relations,
up to some maximum value, which we set at 1800, i.e.\ five times
larger than the previous record $k=360$, set by~\cite{CLZ}.

First, one rules out indices that fail the factorization criterion.
Then, for each surviving $k$, one performs an integer relation search,
at suitably chosen numerical precision, using the constants $\log x$,
$\log\Phi_k(x)$ and a subset of $\{\log\Phi_j(x)\mid j<k\}$
that is consistent with the requirement that
\begin{equation}
N_k=\prod_{0<j<k} N_j^{A_j}
\label{require}
\end{equation}
be satisfiable by rational numbers $A_j$. Moreover, this subset can be further
-- and often greatly -- reduced
by exploiting the cyclotomic relations with indices less than $k$.
In fact, we obtained integer values of $A_j$
that vanish for each $j>d/2$ that does not divide $k$.

For each putative cyclotomic relation thus indicated by a numerically
discovered integer relation between logarithms,
one has merely to use the defining polynomial, to eliminate $x^d$,
via a computer algebra program, hence proving (or if one is very unlucky
disproving) the numerically suggested relation.
This is how we discovered and proved the cyclotomic
relation~(\ref{a630}), with index $k=630$,
in the number field defined by Lehmer's polynomial~(\ref{a10}).
We are also strongly convinced that there is no cyclotomic
relation with $630<k\le1800$ in Lehmer's number field.

Taking the log of a cyclotomic relation, with a real root $x>1$,
one proves the vanishing of a combination of logarithms of the form
\begin{equation}
{\rm Li}_1(x^{-k})-B_0\log x
+\sum_{0<j<k}B_j\,{\rm Li}_1(x^{-j})=0
\label{log}
\end{equation}
where the logarithm ${\rm Li}_1(y):=-\log(1-y)$ is merely the $n=1$ case of
the polylogarithm ${\rm Li}_n(y):=\sum_{r>0}y^r/r^n$, with order $n$,
and the rational numbers $B_j$
are easily obtained from the rational numbers $A_j$
in the cyclotomic relation~(\ref{cycr}).

Next, we define polylogarithmic combinations
\begin{equation}
L_k^{(1)}(n):={\rm Li}_n(x^{-k})/k^{n-1}+B_0(-\log x)^n/n!
+\sum_{0<j<k}B_j\,{\rm Li}_n(x^{-j})/j^{n-1}
\label{plog}
\end{equation}
where the superscript indicates that they all vanish at $n=1$.
In general, these do not vanish at $n=2$.
Rather, one finds~\cite{LL/AMS} that combinations of them
evaluate to rational multiples of $\zeta(2)=\pi^2/6$.

Suppose that several {\bf Q}-linear combinations of the
constructs~(\ref{plog}) evaluate to rational multiples of $\pi^2$
at $n=2$, so that
\begin{equation}
\sum_{k} C^{(2)}_{jk}L_k^{(1)}(2)=D^{(2)}_j\pi^2
\label{Cjk}
\end{equation}
with a rational matrix $C^{(2)}$ yielding a rational vector $D^{(2)}$,
where the superscript indicates that we have exploited
empirical data at order $n=2$.
Then one forms a vector whose components
\begin{equation}
L^{(2)}_j(n)=\sum_{k} C^{(2)}_{jk}L_k^{(1)}(n)
-D^{(2)}_j\frac{\pi^2(-\log x)^{n-2}}{(n-2)!}
\label{v2}
\end{equation}
vanish at $n=2$. At $n=3$, one seeks
{\bf Q}-linear combinations of~(\ref{v2})
that evaluate, empirically, to rational
multiples of $\zeta(3)$. However, these are not yet the constructs
to carry forward to orders $n>3$; one must form combinations
that vanish at $n=3$. In a self-reciprocal number field,
this cannot be done by a subtraction
similar to that in~(\ref{v2}), since $\zeta(3)$ does not appear
in the formula~\cite{LL} for inverting the argument
of a polylogarithm.

Thus the generic iteration is to form combinations
\begin{equation}
L^{(2p+1)}_j(n)=\sum_{k} C^{(2p+1)}_{jk}L^{(2p)}_k(n)
\label{vo}
\end{equation}
that vanish for $n=2p+1$, and then combinations
\begin{equation}
L^{(2p+2)}_j(n)=\sum_{k} C^{(2p+2)}_{jk}L^{(2p+1)}_k(n)
-D^{(2p+2)}_j\frac{\pi^{2p+2}(-\log x)^{n-2p-2}}{(n-2p-2)!}
\label{ve}
\end{equation}
that vanish for $n=2p+2$.

The vital issue is this: how does the number of valid ladders
decrease at each iteration? Don Zagier observed that the answer
depends on the signature of the number field~\cite{DZ}.
Suppose that the polynomial $P(x)$ has $r>0$ real roots
and $s$ pairs of complex roots, so that the degree is
$d=r+2s$. (Note that $r=2$ and $s>0$ in the particular case
of a Salem number.)
Since we restrict attention to reciprocal
polynomials, both $r$ and $d$ are even. (In particular,
the Lehmer polynomial has $r=2$, $s=4$, $d=10$.)
Given $N(2p)$ ladders that evaluate to rational multiples
of $\pi^{2p}$ at order $n=2p$, one expects
\begin{equation}
N(2p+1)=N(2p)+1-d/2
\label{ito}
\end{equation}
ladders that evaluate to rational multiples of $\zeta(2p+1)$
at order $n=2p+1$, and
\begin{equation}
N(2p+2)=N(2p+1)-1-s
\label{ite}
\end{equation}
ladders that evaluate to rational multiples of $\pi^{2p+2}$
at order $n=2p+2$. The $+1$ in~(\ref{ito}) occurs because
one may include the index $k=0$, corresponding to $\zeta(2p+1)$;
the $-1$ in~(\ref{ite}) because this odd zeta value may not
be carried forward. The $-d/2$ in~(\ref{ito}) occurs because
of conditions on the functionally independent real parts of polylogs of
odd order, in a self-reciprocal number field;
the $-s$ in~(\ref{ite}) because of
conditions on the imaginary parts of polylogs of even order.

For a self-reciprocal number field with degree $d=2+2s$
-- and hence for any Salem number -- it follows that
$C$ cyclotomic relations are expected to generate
\begin{equation}
N(n)=C+d/2-\lfloor(d-1)n/2\rfloor
\label{master}
\end{equation}
rational multiples of $\zeta(n)$ at order $n\ge2$.

By way of example, the self-reciprocal number field $\alpha^2-3\alpha+1=0$
has $C=4$ cyclotomic relations, with indices $k=1,6,10,12$.
There are thus $5-\lfloor n/2\rfloor$
valid ladders at orders $n=2\ldots9$. At $n=9$, there is an
integer relation
\begin{equation}
f\,\zeta(9)=
\sum_{j=0}^4g_j\,\pi^{2j}(\log\phi)^{9-2j}
+\sum_{k\in D(\{10,12\})}h_k\,{\rm Li}_{9}(\phi^{-2k})
\label{golden}
\end{equation}
where $\phi:=(1+\sqrt5)/2$ is the golden ratio
and $f$, $g_j$, $h_k$ are essentially unique integers,
with indices $k$ dividing 10 or 12.
In this very simple case, empirical~\cite{LL/AMS} determination of
\begin{equation}
f=2\times3^3\times5\times7^2\times23\times191\times2161
\label{goldf}
\end{equation}
requires only Euclid's algorithm, to find
the rational ratio of each previously
vanishing ladder to the current zeta value.

Much more impressively,
Cohen, Lewin and Zagier~\cite{CLZ} found
71 cyclotomic relations with indices $k\le360$,
for the smallest known Salem number, with $d=10$.
Thus they obtained
$71+5-9\times8=4$ valid ladders at order $n=16$,
yet no relation at $n=17$.

\section{The Lehmer ladder of order 17}

The tables on pages 368--370 of~\cite{LL/AMS}
exhaust the cyclotomic relations with indices $k\le360$
in the Lehmer number field.
It seemed to us peculiarly inconvenient that this tally was precisely
one short of what is needed to generate a rational multiple
of $\zeta(17)$.
It was also clear how to look for ``the one that got away''.
We calculated the norms $N_k$ of $\{\Phi_k(\alpha)\mid360<k\le1800\}$
and found only one candidate with small factors, namely
$N_{630}=N_{126}=5^6$. It seemed likely that
$\Phi_{630}(\alpha)/\Phi_{126}(\alpha)$ was party to a
cyclotomic relation.  Taking logarithms and using PSLQ,
we readily found the numerical relation
\begin{equation}
{\Phi_{630}(\alpha_1)\over\Phi_{126}(\alpha_1)}=
\frac{\alpha_1^{58}-\alpha_1^{55}}{(\alpha_1-1)^5}
\label{simple}
\end{equation}
which was then proven by repeated substitution for $\alpha_1^{10}$.
It entails terms of the form
$\alpha_1^j-1$, where $j$ is one of the 24 divisors of $630$.
These may be halved in number, as in~(\ref{a630}),
by eliminating the 12 divisors
$j\in\{6,7,9,10,18,21,30,42,45,63,70,105\}$, which are
themselves~\cite{LL/AMS} indices of cyclotomic relations.

Proceeding to dilogarithms, we then needed to perform only
6-dimensional searches for integer relations, between the constants
$\zeta(2)$ and $\{L^{(1)}_j(2)|j=6,7,8,9,k\}$,
in the 68 cyclotomic cases with $9<k\le630$.
This resulted in the 67 dilogarithmic ladders of~\cite{LL/AMS}
and one new ladder, namely the integer relation
\begin{eqnarray}
0&=&{\rm Li}_2(\alpha_1^{-630})
-2\,{\rm Li}_2(\alpha_1^{-315})
-3\,{\rm Li}_2(\alpha_1^{-210})
-10\,{\rm Li}_2(\alpha_1^{-126})
-7\,{\rm Li}_2(\alpha_1^{-90})
\nonumber\\
&+&18\,{\rm Li}_2(\alpha_1^{-35})
+84\,{\rm Li}_2(\alpha_1^{-15})
+90\,{\rm Li}_2(\alpha_1^{-14})
-4\,{\rm Li}_2(\alpha_1^{-9})
+339\,{\rm Li}_2(\alpha_1^{-8})
\nonumber\\
&+&45\,{\rm Li}_2(\alpha_1^{-7})
+265\,{\rm Li}_2(\alpha_1^{-6})
-273\,{\rm Li}_2(\alpha_1^{-5})
-678\,{\rm Li}_2(\alpha_1^{-4})
-1016\,{\rm Li}_2(\alpha_1^{-3})
\nonumber\\
&-&744\,{\rm Li}_2(\alpha_1^{-2})
-804\,{\rm Li}_2(\alpha_1^{-1})
-22050\,(\log\alpha_1)^2
+2003\,\zeta(2)
\label{dilog}
\end{eqnarray}
whose index, $k=630$, exceeds anything found previously.
We remark that the coefficients of
$\{{\rm Li}_2(\alpha_1^{-j})\mid 9<j<630\}$ are determined by~(\ref{a630})
and that the empirical coefficient of $\zeta(2)$ is a 4-digit prime,
namely $2003$.

At this juncture, we were faced by a computational dilemma:
how should one process empirical rational data, at orders $n<17$,
so as fastest to determine the final order-17 ladder?
There are two radically opposed strategies:
one systematic, though numerically intensive;
the other interventionist, though requiring less numerical
precision. In the first approach, one takes no heed of the explosion
of primes, such as $2003$ at $n=2$.
Rather, one adopts the simplest procedure of eliminating the predicted number
of indices from the lowest currently available,
as in the case above with $n=2$,
where the indices $k=6,7,8,9$ were eliminated in
passing from 72 cyclotomic relations to 68 dilogarithmic relations.
This already differs from the choice adopted by Lewin in~\cite{LL/AMS},
who chose to eliminate the indices $k=7,8,9,10$, leaving
$k=6$ as a survivor. The latter choice
might, at some intermediate stage, produce integers considerably
smaller than those in our method, yet it is difficult to
automate an objective criterion that will efficiently
limit the growth of scheme-dependent integers at all orders $n<17$.
Since we envisage a unique valid ladder at order 17,
the choice of strategy should not affect the final result.
Rather, it affects the working precision that is required.

Happily, the choice is not crucial,
since state-of-the-art implementation~\cite{ppslq}
of the PSLQ algorithm enables a 6-dimensional search in seconds,
when the relation involves integers with less than 600 digits.
Thus we were able to experiment at lower orders.
The rule of thumb which emerged is this:
sub-optimal intermediate integers, produced by systematic elimination
of the lowest indices currently available,
rarely exceed the squares of those that might be achieved by
laborious optimization.

This suggested that systematic elimination
would be likely to get us fastest to the ultimate goal,
without need for any tuning.
It transpired that the task was indeed as easy as we had supposed,
since the integer $a$ in~(\ref{ais})
has merely 288 digits,
and all the other integers in~(\ref{result}) have less than 300 digits.
At no stage did we encounter, in our systematic approach,
an integer with more than 600 digits, consistent with the rule of thumb.
Thus, in searches with merely 6 constants, we needed less than 4,000-digit
working precision, which placed no significant burden on
MPFUN~\cite{mpfun,mpf90} or PSLQ~\cite{ppslq,cpslq}.

We remark that on a 433 MHz DecAlpha machine it took
139 seconds to compute the constants
$\{{\rm Li}_{17}(\alpha_1^{-k})\mid k\in D({\cal S})\}$ to 4,000 digits.
It then took merely 9 seconds for PSLQ
to find an integer relation between $\zeta(17)$ and the 5
valid ladders that had survived from the previous iteration.
Unravelling all the iterations, we then expressed the final -- and
presumably unique -- result
in the form~(\ref{result}). Of the 117 divisors of~(\ref{DS}), two
are not entailed by cyclotomic relations, namely $k=51$ and $k=53$.
Among the 115 non-zero integers $c_k$, the cyclotomic input
guarantees the triviality of 19 ratios, namely
\begin{eqnarray}
{c_{k}\over c_{2k}}&=&-2^{n-1}\,,\quad k\in
\{41,43,49,59,69,77,91,93,102,106,123,135,143,180,315\}\,;\nonumber\\
{c_{k}\over c_{3k}}&=&-3^{n-1},\quad k\in\{68,210\}\,;\quad
{c_{82}\over c_{246}}=-2\times3^{n-1}\,;\quad
{c_{126}\over c_{630}}=-2\times5^{n-1}
\label{126}
\end{eqnarray}
where, for example, the final ratio results from~(\ref{simple})
and in the particular case of~(\ref{result}) the order is $n=17$.
Moreover, the logarithms of~(\ref {result}) may be removed by replacing
${\rm Li}_n(y)$ by the Rogers-type polylogarithm~\cite{LL/AMS}
\begin{equation}
L_n(y):=\sum_{r=1}^n\left(1-\frac{\delta_{r,1}}{n}\right)
\frac{(-\log|y|)^{n-r}}{(n-r)!}{\rm Li}_r(y)
\label{Rogers}
\end{equation}
which slightly modifies Kummer's~\cite{LL}
$\Lambda_n(-y):=\int_0^y (\log|x|)^{n-1}dx/(1-x)$,
by an inessential normalization
factor, and by a Kronecker delta term, at $r=1$. The latter completes
the process of removal of logs from functional equations, almost
achieved by Kummer. Thus the core of relation~(\ref{result})
is specified by the coefficient~(\ref{ais}) of $\zeta(17)$
and 96 integers $c_k$, from which the full set of 115
is trivially generated. Analyzing these 97 integers in pairs, we found
no pair with a common prime factor greater than $1973$, which divides
both $c_1$ and $c_{57}$. From this, one sees that brute-force
application of PSLQ would require nearly 30,000 digits of precision to
reconstruct the relation, without benefit of further theory. In fact,
4,000 digits were more than enough to ascend the
ladder~(\ref{ito},\ref{ite}), rung by rung.

We presume that the coefficient~(\ref{ais}) of $\zeta(17)$
is essentially unique, making it a remarkable integer in the
theory of polylogarithms. Its 217-digit factor~(\ref{C217})
is certainly composite, yet has so far
resisted factorization\footnote{An ability to factorize the product
of a pair of 100 digit primes might undermine cryptography.}
by Richard Crandall and his colleagues Karl Dilcher,
Richard McIntosh and Alan Powell,
whose large-integer code~\cite{REC} efficiently implements an
elliptic curve method~\cite{RPB}.
Its presence makes it very unlikely that one could discover the final integer
relation~(\ref{result}) with less than 1,000 digits of working precision,
however hard one tried to emulate the feats of~\cite{CLZ}, by interventions
that limit the growth of scheme-dependent integers at lower orders.
By contrast, the ladder given in~\cite{CLZ}, at order $n=16$, involved
integers with no more than 71 digits, which we shall shortly reduce by
11 digits. Nonetheless, we are left speechless with admiration for the
achievement of Henri Cohen in attaining order $n=16$ with only 305-digit
precision, from Pari.

To check~(\ref{result}), we computed more than 59,000 digits of
$\{{\rm Li}_{17}(\alpha_1^{-k})\mid k\in D({\cal S})\}$, using
\begin{equation}
{{\rm Li}_n(\alpha_1^{-k})\over k^n}=\sum_{k|j}{\alpha_1^{-j}\over j^n}
\label{elem}
\end{equation}
which conveniently reduces the 114 cases with $k>1$ to subsets of the
additions for $k=1$. To compute $\zeta(17)$, we set $p=4$ in the identity
\begin{eqnarray}
p\,\zeta(4p+1)&=&\frac{1}{\pi}\sum_{n=0}^{2p+1}
(-1)^n\,(n-\mbox{$\frac12$})\,\zeta(2n)\,\zeta(4p+2-2n)
\nonumber\\&&{}
-2\sum_{n>0}\frac{n^{-4p-1}}{\exp(2\pi n)-1}
\left(p+\frac{\pi n}{1-\exp(-2\pi n)}\right)
\label{z4p}
\end{eqnarray}
which corrects the upper limit of the first sum
and the sign of the second term of the second summand
in Proposition~2 of~\cite{CLZ}.
The corrected companion identity simplifies to
\begin{equation}
\zeta(4p-1)=-\frac{1}{\pi}\sum_{n=0}^{2p}
(-1)^n\,\zeta(2n)\,\zeta(4p-2n)
-2\sum_{n>0}\frac{n^{-4p+1}}{\exp(2\pi n)-1}
\label{z4m}
\end{equation}
At $p=0$, one finds that~(\ref{z4p}) evaluates to $1/4$
and~(\ref{z4m}) to $-1/12$, as expected.

Setting $n=17$ in~(\ref{elem}), we obtained the 115 polylogarithmic constants
to a precision of $3\times2^{16}$ binary digits.
Substituting these in~(\ref{result}), with the integers found by PSLQ
at less than 4,000-digit precision, we reproduced the first 59,157
decimal places of the value of $\zeta(17)$ from~(\ref{z4p}).
The chance of a spurious result is thus less than $10^{-55000}$.

Thanks to PSLQ, we are able to
scrutinize an inference in~\cite{CLZ},
concerning the role of 3617, the numerator
of the Bernoulli number $B_{16}$, and hence of $\zeta(16)/\pi^{16}$.
Appendix~A of~\cite{CLZ} gives one of four empirical
integer relations found at order $n=16$,
with indices $k\le360$. It is of a form
\begin{equation}
\sum_{k\in D({\cal S})} a_k\,L_{16}(\alpha_1^{-k})
=a_0\,\frac{L_{16}(1)}{3617}
\label{3617}
\end{equation}
On the left,
one encounters in~\cite{CLZ} 111 non-zero integers,
with up to 71 digits, and on the right $a_0$ has 75 digits.
This led the authors of~\cite{CLZ} to infer that
ladders do not pick up non-trivial numerators of Bernoulli numbers,
since the relation appears more natural when written in
terms of $\pi^{16}$.
However, such evidence was adduced without knowledge of the
cyclotomic index $k=630$.

Our first remark is that no firm conclusion may be drawn from
Appendix~A of~\cite{CLZ}, since it contains integers with 11 digits
more than is necessary.
We emphasize that this does {\em not\/} indicate
a failure of Pari's LLL~\cite{LLL} algorithm in~\cite{CLZ}.
Rather, it shows that a new analysis is required by~(\ref{a630}).
Using PSLQ, we obtained 5 independent integer relations of the
symmetrical form
\begin{equation}
\sum_{k\in D({\cal S})}s_k\,L_{16}(\alpha_1^{-k}) =
s_0\,L_{16}(1)
\label{no3617}
\end{equation}
with integers $s_k$ having no more than 62 digits and Euclidean norms
in the narrow range
\begin{equation}
2.89\times10^{61}\,<\,\left(s_0^2+\sum_{k>0}s_k^2\right)^{1/2}\,<\,
3.79\times10^{61}
\label{norms}
\end{equation}
We also obtained 5 asymmetrical relations
of the form~(\ref{3617}), with norms
\begin{equation}
2.31\times10^{61}\,<\,\left(a_0^2+\sum_{k>0}a_k^2\right)^{1/2}\,<\,
4.37\times10^{61}
\label{norma}
\end{equation}
Comparison of~(\ref{norms}) with~(\ref{norma}) leaves the issue moot.
Moreover, these ranges were unchanged by application of the LLL algorithm.
An indication of the slight superiority of the asymmetrical
form~(\ref{3617}) is obtained by asking LLL to reduce the restricted
norms that omit the integer $a_0$. This produced 5
asymmetrical relations with
\begin{equation}
0.25\times10^{61}\,<\,\left(\sum_{k>0}a_k^2\right)^{1/2}\,<\,
0.78\times10^{61}
\label{norml}
\end{equation}
The integer relation with the smallest restricted norm is presented in
Table~1, whose 115 integers have no more than 60 digits and
yield the 64-digit integer
\begin{eqnarray}
a_0&=&
2^3\times11\times
1770708910425291120521033962427
\times\nonumber\\&&
23216857398851664164043691705297
\label{a0}
\end{eqnarray}
with factors found after running~\cite{REC} for a total of 466 CPUhours,
on a cluster of 8 machines.

As further indication that powers of $\pi^2$, rather than even zeta values,
are the constants favoured by polylogarithm ladders, we remark that
$b_8$ in~(\ref{result}) is not divisible by 3617, and that $b_6$
is not divisible by 691, the numerator of $\zeta(12)/\pi^{12}$.
Further information on the role of Bernoulli numerators in polylogarithm
ladders could be obtained if one found a base in which there is a unique
valid ladder of order $n=12$ or order $n=16$. We hold out
no lively hope for the latter.
In search of the former, we turn to larger Salem numbers.

\section{Ladders from larger Salem numbers}

There are 47 known~\cite{MJM} Salem numbers less than 1.3.
Of these, 45 exhaust the possibilities
with $\alpha<1.3$ and degree $d<42$.
Of these, merely 6 have degree $d<12$.
Of these 6, we noted that all but one solve equations
of the very simple form
\begin{equation}
x^{4+m}=\frac{Q(1/x)}{Q(x)}
\label{family}
\end{equation}
with $m>0$ and
\begin{equation}
Q(x):=x^3-x-1
\label{Q}
\end{equation}
The case $m=1$ gives Lehmer's number field.
The minimal polynomials of the first five Salem numbers in this family are
\begin{eqnarray}
P_1(\alpha)&=&\alpha^{10}+\alpha^9-\alpha^7-\alpha^6-\alpha^5
-\alpha^4-\alpha^3+\alpha+1\label{P1}\\
P_2(\alpha)&=&\alpha^{10}-\alpha^7-\alpha^5-\alpha^3+1\label{P2}\\
P_3(\alpha)&=&\alpha^{10}-\alpha^8-\alpha^5-\alpha^2+1\label{P3}\\
P_4(\alpha)&=&\alpha^8-\alpha^5-\alpha^4-\alpha^3+1\label{P4}\\
P_5(\alpha)&=&\alpha^{10}-\alpha^8-\alpha^7+\alpha^5-\alpha^3-\alpha^2+1
\label{P5}
\end{eqnarray}
with approximate numerical roots -- and hence Mahler measures -- given by
\begin{eqnarray}
\alpha_1&=&1.1762808182599175065440703384\ldots\label{a1}\\
\alpha_2&=&1.2303914344072247027901779389\ldots\label{a2}\\
\alpha_3&=&1.2612309611371388519466715030\ldots\label{a3}\\
\alpha_4&=&1.2806381562677575967019025327\ldots\label{a4}\\
\alpha_5&=&1.2934859531254541065199098837\ldots\label{a5}
\end{eqnarray}
As $m\to\infty$, one obtains the real root of~(\ref{Q}), namely
\begin{eqnarray}
\alpha_{\infty}&=&
 \left({1+\sqrt{23/27}\over2}\right)^{1/3}
+\left({1-\sqrt{23/27}\over2}\right)^{1/3}
\nonumber\\
&=&1.3247179572447460259609088544\ldots\label{ainf}
\end{eqnarray}
which is the smallest algebraic integer, $\alpha>1$, with conjugates
that all lie inside the unit circle,
i.e.\ the smallest Pisot--Vijayaragharvan (PV) number.

We were alerted to the existence of this family of Salem numbers by
graph theory. Gert Almkvist told the second author of a fascinating remark
on page 247 of the book~\cite{GHJ} by Frederick Goodman, Pierre de la Harpe
and Vaughan Jones, on Coxeter graphs. There they are concerned with the
smallest possible value, greater than 2, for the largest eigenvalue,
$\lambda_{\rm max}$, in the spectrum of a graph.
They classify all finite connected graphs with
$\phi^{1/2}+\phi^{-1/2}\ge\lambda_{\rm max}>2$
where $\phi:=(1+\sqrt5)/2$ is the golden ratio,
and hence prove that the smallest possible
$\lambda_{\rm max}>2$ is obtained from the tree graph ${\rm T}_{2,3,7}$,
formed when three straight lines, with 2, 3, and 7 vertices,
are joined by identifying three univalent vertices.
The characteristic polynomial of its adjacency matrix is
\begin{equation}
P^{\rm T}_{2,3,7}(\lambda)
=\lambda^{10}-9\lambda^{8}+27\lambda^{6}-31\lambda^{4}+12\lambda^{2}-1
\label{T237}
\end{equation}
which is, intriguingly, related to Lehmer's polynomial~(\ref{a10}).
The relation -- which is not correctly stated in~\cite{GHJ} -- is
\begin{equation}
\alpha^5P^{\rm T}_{2,3,7}\left(\alpha^{1/2}+\alpha^{-1/2}\right)
=\alpha^{10}+\alpha^9-\alpha^7-\alpha^6-\alpha^5-\alpha^4-\alpha^3+\alpha+1
\label{T237a}
\end{equation}
which proves that the largest eigenvalue of ${\rm T}_{2,3,7}$ is
\begin{equation}
\lambda_1^{}:=\alpha_1^{1/2}+\alpha_1^{-1/2}=
2.0065936183460167326505159176\ldots
\label{la}
\end{equation}
Thanks to help from Gert Almkvist,
we proved that this generalizes to the relation
\begin{equation}
\lambda_m^{}=\alpha_m^{1/2}+\alpha_m^{-1/2}
\label{lag}
\end{equation}
between the largest eigenvalue of the tree graph ${\rm T}_{2,3,6+m}$
and the Salem number obtained from~(\ref{family}), at any $m>0$.

Thus this family of graphs generates a monotonically increasing
family of Salem numbers, starting at the smallest yet known
and tending to the provably smallest PV number~(\ref{ainf}).
It is proven~\cite{GHJ}
that~(\ref{la}) is the smallest value, $\lambda_{\rm max}>2$,
for the largest eigenvalue of any graph. Sadly,
this does not prove that Lehmer's $\alpha_1$ is the smallest
Salem number.

Since Lehmer's number field
generates the new record, $n=17$,
for the order of a polylogarithmic ladder, we thought it interesting
to examine ladders based on further members of the family of graphs
${\rm T}_{2,3,6+m}$.
By way of benchmark, we recall
that the square of the golden ratio
generates ladder~(\ref{golden}), with order $n=9$ and index $k=12$.
The smallest PV number~(\ref{ainf}) also reaches $n=9$,
with index $k=42$~\cite{LL/AMS}.
In~\cite{z3z5}, order $n=11$ was attained, from cyclotomic
relations that also led to algorithms for finding
the ten millionth hexadecimal digits of $\zeta(3)$ and $\zeta(5)$,
without computing previous digits.
To our knowledge, $n=11$ had been bettered, heretofore, only
in the Lehmer number field.

We were intrigued to know the maximum orders achievable
in bases derived from trees ${\rm T}_{2,3,6+m}$ with $m>1$.
Accordingly, we studied the cyclotomic relations of~(\ref{family}),
with $m<6$, and found $N_m(n)$ ladders, for
base $\alpha_m$ at order $n>1$, with
\begin{eqnarray}
N_1(n)&=&77-\lfloor9n/2\rfloor\quad(n\le17)\label{n1}\\
N_2(n)&=&54-\lfloor9n/2\rfloor\quad(n\le11)\label{n2}\\
N_3(n)&=&49-\lfloor9n/2\rfloor\quad(n\le10)\label{n3}\\
N_4(n)&=&44-\lfloor7n/2\rfloor\quad(n\le12)\label{n4}\\
N_5(n)&=&43-\lfloor9n/2\rfloor\quad(n\le9)\label{n5}
\end{eqnarray}
The number of
cyclotomic relations decreases monotonically, as the Mahler measure
increases from $\alpha_1$ to $\alpha_5$. However $\alpha_4$
carries the advantage of having the lowest degree, $d=8$, among
Salem numbers from trees ${\rm T}_{2,3,6+m}$.
The explanation is that the characteristic polynomial
of the adjacency matrix of tree ${\rm T}_{2,3,10}$ is
\begin{equation}
P^{\rm T}_{2,3,10}(\lambda)=(\lambda^4-4\lambda^2+2)
P^{\rm T}_{2,4,5}(\lambda)
\label{T2310}
\end{equation}
with the same largest
eigenvalue as for tree ${\rm T}_{2,4,5}$, which is
the first of a family of trees ${\rm T}_{2,4,4+m}$
whose largest eigenvalues
are $\beta_m^{1/2}+\beta_m^{-1/2}$, where the Salem
number $\beta_m$ solves
\begin{equation}
x^{2+m}=\frac{R(1/x)}{R(x)}
\label{other}
\end{equation}
with $m>0$ and
\begin{equation}
R(x):=x^3-x^2-1
\label{R}
\end{equation}
These Salem numbers increase monotonically, from $\beta_1=\alpha_4$, to
\begin{eqnarray}
\beta_{\infty}&=&{1\over
\left({\sqrt{31/27}+1\over2}\right)^{1/3}
-\left({\sqrt{31/27}-1\over2}\right)^{1/3}}
\nonumber\\
&=&1.4655712318767680266567312252\ldots
\label{binf}
\end{eqnarray}
which is the PV number obtained from~(\ref{R}).
As seen from~(\ref{n4}), the first of these Salem numbers gives
valid ladders up to $n=12$. These entail 62 indices that divide elements of
\begin{eqnarray}
&&{\cal T}:=
\{23, 24, 33, 40, 45, 50, 54, 55, 56, 60, 62,\nonumber\\&&
  64, 68, 75, 78, 84, 88,102,105,114,140,252\}
\label{DT}
\end{eqnarray}
Finding two valid ladders at order $n=12$,
we obtained no definitive answer regarding the appearance of 691,
the numerator of $\zeta(12)/\pi^{12}$.

This study of graphical number fields led to the observation that
the Salem number
\begin{equation}
\alpha_{\rm not}=1.216391661138265091626806311199463327722253606570\ldots
\label{not}
\end{equation}
that solves
\begin{equation}
\alpha^{10}-\alpha^6-\alpha^5-\alpha^4+1=0
\label{notpol}
\end{equation}
is rather special. It is the sole Salem
number $\alpha<1.3$ with degree $d<12$ for which $\alpha^{1/2}+\alpha^{-1/2}$
is {\em not\/} the largest eigenvalue of the adjacency matrix of a graph.
We note that
Gary Ray~\cite{GAR} identified~(\ref{notpol}) as a potentially
fruitful source of polylogarithm identities, along
with~(\ref{P1},\ref{P4}).

Accordingly, we sought
its cyclotomic relations and hence its polylogarithmic ladders, finding
$59-\lfloor9n/2\rfloor$ valid ladders at orders $n=2\ldots13$. At $n=13$,
we found an integer relation of the form
\begin{equation}
u\,\zeta(13)=
\sum_{j=0}^6v_j\,\pi^{2j}(\log\alpha_{\rm not})^{13-2j}
+\sum_{k\in D({\cal U})}w_k\,{\rm Li}_{13}(\alpha_{\rm not}^{-k})
\label{lad13}
\end{equation}
with integers (see footnote~1)
$u$, $v_j$ and $w_k$, where $k$ runs over 86 divisors
of elements of
\begin{eqnarray}
&&{\cal U}:=
\{32,38,40,43,48,50,54,56,57,58,60,72,75,84,88,90,
\nonumber\\&&\phantom{90,}
92,110,118,124,126,132,136,156,170,204,210,234,294\}
\label{DU}
\end{eqnarray}
The coefficient of $\zeta(13)$ was factorized as follows
\begin{eqnarray}
u&=&2^3\times3^5\times5^2\times7^2\times11\times13\times19
\times43\times107\times1789\times3413\nonumber\\&&{}
\times350162215794091
\times31692786317928349
\times P_{94}
\label{ris}
\end{eqnarray}
where the 94-digit factor
\begin{eqnarray}
&&P_{94}:=
80131278764880863222650485358197334145287799\nonumber\\&&
30835715165005435238201238560738060500528870828871
\label{P94}
\end{eqnarray}
is (very) probably prime.

We believe that~(\ref{lad13}),
with index $k=294$, is the unique valid ladder of order $n=13$
in the number field~(\ref{notpol}). In our investigations,
the only number field that yielded orders $n>13$ was Lehmer's.
Thanks to the index $k=630$ in~(\ref{a630}), the latter attains
order $n=17$. We expect this record to abide.

\subsubsection*{Acknowledgements}
DJB thanks Gert Almkvist, for an enjoyable stay in Lund,
where the relation to graph theory emerged;
Richard Crandall, for advice on factorization;
Petr Lisonek, for discussions at Simon Fraser University, from which
the idea of using PSLQ on the Lehmer number field originated;
Chris Wigglesworth for expert management of a cluster of Dec\-Alpha
machines at the Open University; and Don Zagier, for discussions in Bonn
and Vienna, which helpfully emphasized the reliability of~(\ref{master}).

\newpage
\begin{center}{\bf Table~1}:
The 115 integers $a_k$ in~(\ref{3617}) yield~(\ref{a0}) and\\
have 11 digits less than those in Appendix~A of~\cite{CLZ}.
\end{center}
$$\begin{array}{|r|r|}\hline k&a_k\hspace{7cm}{}\\\hline
1&657517430619563136979927560311907031673621211923757556039680\\
2&359886000860172792447593380451304667509342919316395916742210\\
3&-141684894142517938914234443949080358856724456088731117486080\\
4&196781701168242866806602275650062000534697577654107935879825\\
5&745935348767332443069712397147123742304948147537862045466624\\
6&405488816708650056442253795413843043306948102060187723235830\\
7&-416406976403628958702584514595218236213894104883105761525760\\
8&-547505077081618409188348198016100365604103632542087633890050\\
9&-168975342608695257899656804048229345162975246342759199539200\\
10&-230249343614447384417237255880683163036311216601594850723110\\
11&-686361499027198737843324163158516540202013109745360933355520\\
12&-555211669354532402241841641332915138353295668921519456494790\\
13&192815301162191732499187837598868261035304857047105259765760\\
14&679734532959241358849438888514522682758692394636208781152890\\
15&542262843988600709919208184100477463027334034487684580835328\\
16&424235600610741544849445172143286300309734125137633061204035\\
17&523715028246745518068354598319413022191932310750064994549760\\
18&391765228725624091642051568538398618159246632370196009838220\\
19&-198765352117361763621609618366336349702431673161386991943680\\
20&-681771795414951811771620988305943452975785965232369963176182\\
21&-711962908052232232505015799802945837905471681104673666826240\\
22&-474711713652489086239653502137393573002849538293395396856580\\
23&-305563488793803341514102820998121597667241109957343639961600\\
24&-72481617150888165437226250731733584917169373028647243993490\\
25&277441845925552929689261105309806909812357891559963458994176\\
26&560813990488672759149774541727869507575569280403768300321140\\
27&529382955893422715273603248526440343418480284327032762204160\\
28&320081042366233059912738009855303427011191818027491322896270\\
29&83494538290613239874707297912762677581023131627678971658240\\
30&-98012487241465589885244822793749074739104108244040188504014\\
31&-259850303587986534656435580762491706946152146552951829626880\\
32&-312381473762254154637700526528599595535218857011814459284405\\
33&-246779089811622911261098199061852288809784802679547056947200\\
34&-101370164345933351431594273523846287706927174737460926584690\\
35&10117727926311381644405806386440429351025303665206665019392\\
36&92813012460581905456644049699805057100444523762661712331760\\
37&126324305118118269024517948173064327264644326236420977131520\\
38&115549700605836091570299617434612530310676419462990155960060\\
39&55524385597088576211029994800590286128713854655226225950720
\end{array}$$

\newpage
$$\begin{array}{|r|r|}
40&8681384820081409297328583655855503616180543689183078659790\\
41&-19635305068579015455281545202962471423318427323980082053120\\
42&-25853290618480167202972487269611444701702175335712826764200\\
43&-31825441142023707248687536800399542265330919953456915087360\\
44&-21238990403232997150237193584099678135190886191312613084840\\
45&-5941171482587742832139996778309010890215227722504605532160\\
46&6867449991711733361905911701040369140850427873891452167820\\
47&2683941033568052348250476195016210444676454909493450178560\\
48&750276335951027139956797045951168533453325950288012793755\\
49&902953609381753090186332205561649133432166531686386565120\\
50&2845791617893774854256316723847280908911507508173997476472\\
52&-211518744564361609748949526157272452992328540189820187580\\
54&1586305520716919558891270499363536151171957733664553920190\\
55&317399573539556012199105353632753466414027589073604640768\\
56&-411955065839184813654206085597826201289284276111276432380\\
57&-633034566415505822568457398619791655535556228581668945920\\
59&-422235192545960711449176250527328955777038182764768133120\\
60&-100775809060677875668353632672084730591567385836002358340\\
62&103161260941635325589940465107263862028544860193511288430\\
63&-71239824525936132532239108327787334379823490760552611840\\
64&-758961502664172364046541735177189723413489520451845090\\
65&-31614329666390498308045970054858836918855503060815904768\\
66&82590308807941016435157163298430216384039435573691923590\\
68&-6111893964221453989555021260745192324274991457297917375\\
69&36723207357049059339848864763934042997356103557409013760\\
70&56587501615035807034718738004427311071764030581924717666\\
72&6246571495623829969272301654360557650197045971497886770\\
74&4662327569727635575068647456890293743838455770712101950\\
75&14887216546252872098059151093987262212120215877194153984\\
76&4933420845476441755128409798981902904608462970609785585\\
77&-3718375935394449285719242013677644130212998760206499840\\
78&797543027544694301981107935492458580524804545433769960\\
80&-756655048589380952520516473920663602783235627959554055\\
82&599221956438568586892136999602126203104200052611696840\\
84&1711070904277206791256558392136135387624754115179524010\\
86&971235386414297706563950708020005562296475828657742770\\
90&712743581327703411348920012655781994420044560878114988\\
91&581533337739120206692138569579040767870857035407360000
\end{array}$$

\newpage
$$\begin{array}{|r|r|}
92&281171208042281370705030709443218746311661931924480940\\
93&-123474844527238243730367824366624184577805521172234240\\
96&45358670699766680263040092382358986091316360899268450\\
98&-27555957317558382879221563890431186933354691518749590\\
102&-13957477138823786670980509990907214192819210555392000\\
105&-51592349989081948844625156379637404692389576604581888\\
106&-49757317433194357135489925493889989029913561751388160\\
108&4104159613575523414703444177171856613076465917080545\\
110&436503492818404043948452074062000221710038770175550\\
118&12885595475645773664830818192362333855500432823631840\\
120&-4458271535866469116508257194072443804522804186654722\\
123&684209085353296089217162854388925631175894697902080\\
124&-2568170165780877908485719228826518196430780458722240\\
126&-1238638691691188975878629542803876608020141601562500\\
130&1661804350411188506429822540247968366095406030963502\\
132&-1633602118796242535519471886453758168320481063983580\\
135&1847684862191853167535965349249779010200990442520576\\
138&-1120703349519319437861598656125916839518924058758820\\
143&174515941261835218319905200778076418833749415362560\\
144&-317934100357607293764664092299747258866505243251140\\
154&113475828106520058768287414968189823309722862555130\\
160&-10275921042899706253486332862758642766530474426501\\
165&-23116946816383260126447835448067026267507493044224\\
175&-32347514260620927078585304389860102844649182658560\\
180&-4479106615252570847623887867187960984130837217280\\
182&-17746989066745611776493486620454124996058869488750\\
186&3768153214332221793529291515094732195367600133430\\
195&-2675778372359720766614632517822289310362850099200\\
204&425948399011956380339981383999853948755469072125\\
210&-291194657074029253989297333815629327367478738414\\
212&1518472822057933262191465011410216950375780082745\\
240&76199283284999932540837893030287171626166802389\\
246&-20880404216104006628941737499662037084225302060\\
270&-56386867132319737778807536293023041082793897782\\
286&-5325803871515967355954138207338757898979169170\\
315&-664989084046735447844305843955260271683239936\\
360&136691486061174647449459468603148223392664710\\
630&20293856324668440180795466429298714345802\\\hline
\end{array}$$

\raggedright

\end{document}